\documentclass[11pt]{amsart}
\usepackage{mathrsfs}
\usepackage{amssymb,latexsym}
\usepackage{pstcol,pstricks,color}

 \numberwithin{equation}{section}

\newtheorem{theorem}{Theorem}[section]

\newtheorem{corollary}[theorem]{Corollary}

\def\qed{\hfill $\Box$}
\def\pf{\noindent {\it Proof.} }

\title{ Four transformations on the Catalan triangle }

\begin{document}
\maketitle
\begin{center}
Yidong Sun$^\dag$ and Fei Ma

%%%%\thanks{Corresponding author: Yidong Sun, sydmath@yahoo.com.cn.}%%

Department of Mathematics, Dalian Maritime University, 116026 Dalian, P.R. China\\[5pt]

{\it  Emails: $^\dag$sydmath@dlmu.edu.cn }

\end{center}\vskip0.2cm

\subsection*{Abstract} In this paper, we define four transformations on the classical Catalan
triangle $\mathcal{C}=(C_{n,k})_{n\geq k\geq 0}$ with $C_{n,k}=\frac{k+1}{n+1}\binom{2n-k}{n}$.
The first three ones are based on the determinant and the forth is utilizing the permanent of a square matrix.
It not only produces many known and new identities
involving Catalan numbers, but also provides a new viewpoint on combinatorial triangles.

\medskip

{\bf Keywords}:  Catalan number, ballot number, Catalan triangle, determinant, permanent.

\noindent {\sc 2000 Mathematics Subject Classification}: Primary
05A19; Secondary 05A15, 15A15.

{\bf \section{Introduction} }

In 1976, by a nice interpretation in terms of pairs of paths on a lattice $\mathbb{Z}^2$,
Shapiro \cite{ShapA} first introduced the Catalan triangle $\mathcal{B}=(B_{n,k})_{n\geq k\geq 0}$ such that
$B_{n,k}=\frac{k+1}{n+1}\binom{2n+2}{n-k}$ and obtained
\begin{eqnarray}
\sum_{k=0}^{n}B_{n,k}\hskip-.22cm &=&\hskip-.22cm (2n+1)C_{n}, \nonumber  \\
\sum_{k=0}^{\min\{m,n\}}B_{n,k}B_{m,k}\hskip-.22cm &=&\hskip-.22cm C_{m+n+1}, \label{eqn 1.1.1}
\end{eqnarray}
where $C_n=\frac{1}{n+1}\binom{2n}{n}=\frac{1}{2n+1}\binom{2n+1}{n}$ is the $n$th Catalan number.
Table 1 illustrates this triangle for small $n$ and $k$ up to $6$. Note that the entries in the first column of the Catalan triangle
$\mathcal{B}$ are indeed the Catalan numbers $B_{n,0}=C_{n+1}$, which is the reason why $\mathcal{B}$ is called the Catalan triangle.

\begin{center}
\begin{eqnarray*}
\begin{array}{c|cccccccc}\hline
n/k & 0   & 1    & 2    & 3    & 4    & 5     & 6   \\\hline
  0 & 1   &      &      &      &      &       &     \\
  1 & 2   & 1    &      &      &      &       &      \\
  2 & 5   & 4    & 1    &      &      &       &      \\
  3 & 14  & 14   & 6    &  1   &      &       &      \\
  4 & 42  & 48   & 27   &  8   &  1   &       &     \\
  5 & 132 & 165  & 110  &  44  &  10  &  1    &     \\
  6 & 429 & 572  & 429  &  208 &  65  &  12   & 1   \\\hline
\end{array}
\end{eqnarray*}
Table 1. The values of $B_{n,k}$ for $n$ and $k$ up to $6$.
\end{center}

Since then, much attentions has been paid to the Catalan triangle and its generalizations. In 1979, Eplett \cite{Eplett} deduced
the alternating sum in the $n$th row of $\mathcal{B}$, namely,
\begin{eqnarray*}
\sum_{k=0}^{n}(-1)^{k}B_{n,k}\hskip-.22cm &=&\hskip-.22cm  C_{n}.
\end{eqnarray*}
In 1981, Rogers \cite{RogersA} proved that a generalization of Eplett's identity holds in any renewal array.
In 2008, Guti¨¦rrez et al. \cite{Gutierrez} established three summation identities
%\begin{eqnarray*}
%\sum_{k=0}^{n}(k+1)B_{n,k}^2 \hskip-.22cm &=&\hskip-.22cm  \binom{n+3}{2}C_{n+1}C_{n}, \\
%\sum_{k=0}^{n}(k+1)^2B_{n,k}^2 \hskip-.22cm &=&\hskip-.22cm  (3n+1)C_{2n}, \\
%\sum_{k=0}^{i}(n+2k+2-i)B_{n,k}B_{n,n+k-i} \hskip-.22cm &=&\hskip-.22cm  (n+2)\binom{2n}{i}C_{n+1}, \ (0\leq i\leq n),
%\end{eqnarray*}
and proposed as one of the open problems to evaluate the moments $\Omega_m=\sum_{k=0}^{n}(k+1)^{m}B_{n,k}^2$.
Using the WZ-theory (see \cite{Petkov, Zeilberger}), Miana and Romero computed $\Omega_m$ for $1\leq m\leq 7$.
%and deduced that
%\begin{eqnarray*}
%\sum_{k=0}^{i}(n+2k+2-i)^3B_{n,k}B_{n,n+k-i} \hskip-.22cm &=&\hskip-.22cm  (n+2)((n-i)^2+4(n+1))\binom{2n}{i}C_{n+1}, \ (0\leq i\leq n).
%\end{eqnarray*}
Later, based on the symmetric functions and inverse series relations with combinatorial computations, Chen and Chu \cite{ChenChu}
resolved this problem in general. By using the Newton interpolation formula, Guo and Zeng \cite{GuoZeng} generalized the recent
identities on the Catalan triangle $\mathcal{B}$ obtained by Miana and Romero \cite{Miana} as well as those of Chen and Chu \cite{ChenChu}.

Some alternating sum identities on the Catalan triangle $\mathcal{B}$ were established by Zhang and Pang \cite{ZhangPang}, they
showed that the Catalan triangle $\mathcal{B}$ can be factorized as the product of the Fibonacci matrix and
a lower triangular matrix, which makes them build close connections among $C_n$, $B_{n,k}$ and the Fibonacci numbers.
Motivated by a matrix identity related to the Catalan triangle  $\mathcal{B}$ \cite{ShapGet},
Chen et al. \cite{ChenLi} derived many nice matrix identities on weighted partial Motzkin paths.
Aigner \cite{AignerA}, in another direction, came up with the admissible matrix, a kind of generalized Catalan triangle,
and discussed its basic properties. The numbers in the first column of the admissible matrix are called Catalan-like numbers, which are
investigated in \cite{AignerC} from combinatorial views.

The admissible matrix $\mathcal{A}=(A_{n,k})_{n\geq k\geq 0}$ associated to the Catalan triangle  $\mathcal{B}$ is defined by
$A_{n,k}=\frac{2k+1}{2n+1}\binom{2n+1}{n-k}$, which is considered by Miana and Romero \cite{MianaRom} by
evaluating the moments $\Phi_m=\sum_{k=0}^{n}(2k+1)^{m}A_{n,k}^2$. Table 2
illustrates this triangle for small $n$ and $k$ up to $6$.
\begin{center}
\begin{eqnarray*}
\begin{array}{c|cccccccc}\hline
n/k & 0   & 1    & 2    & 3    & 4    & 5     & 6     \\\hline
  0 & 1   &      &      &      &      &       &      \\
  1 & 1   & 1    &      &      &      &       &      \\
  2 & 2   & 3    & 1    &      &      &       &       \\
  3 & 5   & 9    & 5    &  1   &      &       &       \\
  4 & 14  & 28   & 20   &  7   &  1   &       &       \\
  5 & 42  & 90   & 75   &  35  &  9   &  1    &      \\
  6 & 132 & 297  & 275  &  154 &  54  &  11   & 1    \\ \hline
\end{array}
\end{eqnarray*}
Table 2. The values of $A_{n,k}$ for $n$ and $k$ up to $6$.
\end{center}

The mergence of the two triangles $\mathcal{A}$ and $\mathcal{B}$ forms the third triangle $\mathcal{C}=(C_{n,k})_{n\geq k\geq 0}$, defined by
the ballot numbers
$$C_{n,k}=\frac{k+1}{2n-k+1}\binom{2n-k+1}{n-k}=\frac{k+1}{n+1}\binom{2n-k}{n}.$$
The triangle $\mathcal{C}$ is also called the ``Catalan triangle" in the literature, despite it has
the most-standing form $\mathcal{C}'=(C_{n,n-k})_{n\geq k\geq 0}$, see for examples \cite{BarcVerri, Forder, KitLies, MianaRom, Sloane}.
Table 3 illustrates this triangle for small $n$ and $k$ up to $7$.
\begin{center}
\begin{eqnarray*}
\begin{array}{c|cccccccc}\hline
n/k & 0   & 1   & 2    & 3    & 4    & 5    & 6     & 7   \\\hline
  0 & 1   &     &      &      &      &      &       &     \\
  1 & 1   & 1   &      &      &      &      &       &     \\
  2 & 2   & 2   & 1    &      &      &      &       &      \\
  3 & 5   & 5   & 3    & 1    &      &      &       &      \\
  4 & 14  & 14  & 9    & 4    &  1   &      &       &      \\
  5 & 42  & 42  & 28   & 14   &  5   &  1   &       &     \\
  6 & 132 & 132 & 90   & 48   &  20  &  6   &  1    &     \\
  7 & 429 & 429 & 297  & 165  &  75  &  27  &  7    & 1   \\\hline
\end{array}
\end{eqnarray*}
Table 3. The values of $C_{n,k}$ for $n$ and $k$ up to $7$.
\end{center}

Clearly,
$$A_{n,k}=C_{n+k, 2k}\ \mbox{and} \  B_{n,k}=C_{n+k+1, 2k+1}.$$
Three relations, $C_{n,0}=C_n$, $C_{n+1,1}=C_{n+1}$ and $\sum_{k=0}^{n}C_{n,k}=C_n$
bring the Catalan numbers and the ballot numbers in correlation \cite{AignerB, Hilton, RogersB}. Many properties of the Catalan numbers can be
generalized easily to the ballot numbers, which have been studied intensively by Gessel \cite{Gessel}. The combinatorial interpretations
of the ballot numbers can be found in \cite{AignerC, ChenLi, Comtet, EuLiuYeh, Knuth, MerliSprug, ShapA, SunMa}.
It was shown by Ma \cite{ShiMeiMa} that the Catalan's triangle $\mathcal{C}$ can be generated by context-free grammars in three variables.

The Catalan triangles $\mathcal{B}$ and $\mathcal{C}$ often arise as examples of the infinite matrix associated
to a generating tree \cite{BacchMerli, Merlini, MerlVer}. In the theory of Riordan arrays \cite{ShapB, ShapGet, Sprug},
much interest has been taken in the three
triangles $\mathcal{A}, \mathcal{B}$ and $\mathcal{C}$, see
\cite{ChenLi, CheonJin, CheonKim, CheonKimShap, Luzona, Merlini, Sprugnoli, SunMa}. In fact,
$\mathcal{A}, \mathcal{B}$ and $\mathcal{C}$ are Riordan arrays
\begin{eqnarray*}
\mathcal{A}=(C(t), tC(t)^2), \ \ \mathcal{B}=(C(t)^2, tC(t)^2),\ \mbox{and}\  \mathcal{C}=(C(t), tC(t)),
\end{eqnarray*}
where $C(t)=\frac{1-\sqrt{1-4t}}{2t}$ is the generating function for the Catalan numbers $C_n$.

Recently, Sun and Ma \cite{SunMa} studied the sums of minors of second order of
$\mathcal{M}=(M_{n,k}(x,y))_{n\geq k\geq 0}$, a class of infinite lower triangles
related to weighted partial Motzkin paths and obtained the main theorem.

\begin{theorem}
For any integers $n, r\geq 0$ and $m\geq \ell\geq 0$, set $N_r=\min\{n+r+1, m+r-\ell\}$.
Then there hold
\begin{eqnarray}
\sum_{k=0}^{N_r}\det\left(\begin{array}{cc}
M_{n,k}(x,y)   & M_{m,k+\ell+1}(x,y) \\[5pt]
M_{n+r+1,k}(x,y) & M_{m+r+1,k+\ell+1}(x,y)
\end{array}\right)
\hskip-.22cm &=&\hskip-.22cm \sum_{i=0}^{r} M_{n+i,0}(x,y)M_{m+r-i,\ell}(y,y). \label{eqn 1.1.2}
\end{eqnarray}
\end{theorem}

Recall that a {\em partial Motzkin path} is a lattice path from $(0, 0)$ to
$(n, k)$ in the $XOY$-plane that does not go below the $X$-axis and
consists of up steps $\mathbf{u}=(1, 1)$, down steps $\mathbf{d}=(1, -1)$ and
horizontal steps $\mathbf{h}=(1, 0)$. A {\em weighted partial Motzkin path} (not the same as in \cite{ChenLi})
is a partial Motzkin path with the weight assignment that the all up
steps and down steps are weighted by $1$, the horizontal steps are
endowed with a weight $x$ if they are lying on $X$-axis, and endowed
with a weight $y$ if they are not lying on $X$-axis. The {\em
weight} $w(P)$ of a path $P$ is the product of the weight of all its
steps. The {\em weight} of a set of paths is the sum of the total
weights of all the paths. Denote by $M_{n,k}(x,y)$ the weight sum of the set $\mathcal{M}_{n,k}(x,y)$ of
all weighted partial Motzkin paths ending at $(n,k)$.
The triangle $\mathcal{M}$ can reduce to $\mathcal{A}, \mathcal{B}$ and $\mathcal{C}$ when the parameters $(x,y)$ are specalized, namely,
\begin{eqnarray*}
A_{n,k}=M_{n,k}(1,2),\ B_{n,k}=M_{n,k}(2,2)  \ \mbox{and} \  C_{n,k}=M_{2n-k, k}(0,0).
\end{eqnarray*}
In this paper, we define four transformations on the Catalan triangle $\mathcal{C}$, and derive some known and new identities involving
Catalan numbers. The first three transformations are special cases derived from (\ref{eqn 1.1.2}) which are presented in Section 2 and 3
respectively. Section 4 is devoted to the forth one defined by the
permanent of a square matrix, and gives a general result on the triangle $\mathcal{M}$.

\vskip.5cm
\section{The first two transformation on the Catalan triangle}

Let $\mathcal{X}=(X_{n,k})_{n\geq k\geq 0}$ and $\mathcal{Y}=(Y_{n,k})_{n\geq k\geq 0}$
be the infinite lower triangles defined on the Catalan triangle $\mathcal{C}$ respectively by
\begin{eqnarray*}
X_{n,k} \hskip-.22cm &=&\hskip-.22cm  \rm{det}\left(\begin{array}{cc}
C_{n+k, 2k}    & C_{n+k,2k+1}  \\[5pt]
C_{n+k+1, 2k}  & C_{n+k+1,2k+1}
\end{array}\right), \\
Y_{n,k} \hskip-.22cm &=&\hskip-.22cm  \rm{det}\left(\begin{array}{cc}
C_{n+k+1, 2k+1}    & C_{n+k+1,2k+2}  \\[5pt]
C_{n+k+2, 2k+1}  & C_{n+k+2,2k+2}
\end{array}\right).
\end{eqnarray*}
Table 2.1 and 2.2 illustrate these two triangles $\mathcal{X}$ and $\mathcal{Y}$ for small $n$ and $k$ up to $5$, together
with the row sums. It indicates that the two transformations contact the row sums of $\mathcal{X}$ and $\mathcal{Y}$
with the first two columns of $\mathcal{C}$.
\begin{center}
\begin{eqnarray*}
\begin{array}{c|cccccc|c}\hline
n/k & 0   & 1   & 2    & 3    & 4    & 5    & row\ sums           \\\hline
  0 & 1   &     &      &      &      &      &  1=1^2              \\
  1 & 0   & 1   &      &      &      &      &  1=1^2             \\
  2 & 0   & 3   & 1    &      &      &      &  4=2^2             \\
  3 & 0   & 14  & 10   & 1    &      &      &  25=5^2            \\
  4 & 0   & 84  & 90   & 21   &  1   &      &  196=14^2            \\
  5 & 0   & 594 & 825  & 308  & 36   &  1   &  1764=42^2         \\\hline
\end{array}
\end{eqnarray*}
Table 2.1. The values of $X_{n,k}$ for $n$ and $k$ up to $5$, together with the row sums.
\end{center}
\begin{center}
\begin{eqnarray*}
\begin{array}{c|cccccc|c}\hline
n/k & 0   & 1   & 2    & 3    & 4    & 5      & row\ sums     \\\hline
  0 & 1   &     &      &      &      &        &  1=1\times 1            \\
  1 & 1   & 1   &      &      &      &        &  2=1\times 2          \\
  2 & 3   & 6   & 1    &      &      &        &  10=2\times 5          \\
  3 & 14  & 40  & 15   & 1    &      &        &  70=5\times 14         \\
  4 & 84  & 300 & 175  & 28   &  1   &        &  588=14\times 42         \\
  5 & 594 & 2475&1925  &504   & 45   &  1     &  5544=42\times 132      \\\hline
\end{array}
\end{eqnarray*}
Table 2.2. The values of $Y_{n,k}$ for $n$ and $k$ up to $5$, together with the row sums.
\end{center}

More generally, we obtain the first result which is a consequence of Theorem 1.1.

\begin{theorem} For any integers $m\geq \ell\geq 0$ and $n\geq 0$, set $N=\min\{n+1, m-\ell\}$, then there hold
\begin{eqnarray}
\sum_{k=0}^{N}\det\left(\begin{array}{cc}
C_{n+k,2k}      & C_{m+k,2k+\ell+1}  \\[5pt]
C_{n+k+1,2k}    & C_{m+k+1,2k+\ell+1}
\end{array}\right)
\hskip-.22cm &=&\hskip-.22cm C_{n}C_{m,\ell}, \label{eqn 2.1.1} \\
\sum_{k=0}^{N}\det\left(\begin{array}{cc}
C_{n+k+1,2k+1}      & C_{m+k+1,2k+\ell+2}  \\[5pt]
C_{n+k+2,2k+1}      & C_{m+k+2,2k+\ell+2}
\end{array}\right)
\hskip-.22cm &=&\hskip-.22cm C_{n+1}C_{m,\ell}, \label{eqn 2.1.2}
\end{eqnarray}
or equivalently,
\begin{eqnarray}
C_{m,\ell}C_{n} \hskip-.22cm &=&\hskip-.22cm
\sum_{k=0}^{N}\frac{ (2k+1)(2k+\ell+2)\lambda_{n,k}(m,\ell)}{(2n+1)_3(2m-\ell)_3}\binom{2n+3}{n-k+1}\binom{2m-\ell+2}{m-k-\ell}, \label{eqn 2.1.3}\\
C_{m,\ell}C_{n+1} \hskip-.22cm &=&\hskip-.22cm
\sum_{k=0}^{N}\frac{ (2k+2)(2k+\ell+3)\mu_{n,k}(m,\ell)}{(2n+2)_3(2m-\ell+1)_3}\binom{2n+4}{n-k+1}\binom{2m-\ell+3}{m-k-\ell}, \label{eqn 2.1.4}
\end{eqnarray}
where $\lambda_{n,k}(m,\ell)=(2m-\ell)(2m-\ell+1)(n-k+1)(n+k+2)-(2n+1)(2n+2)(m-\ell-k)(m+k+2)$,
$\mu_{n,k}(m,\ell)=(2m-\ell+1)(2m-\ell+2)(n-k+1)(n+k+3)-(2n+2)(2n+3)(m-\ell-k)(m+k+3)$ and  $(x)_k=x(x+1)\cdots(x+k-1)$ for $k\geq 1$
with $(x)_0=1$.
\end{theorem}
\pf Setting $(x,y)=(0,0)$ and $r=1$, replacing $n,m$ respectively by $2n, 2m-\ell-1$ in (\ref{eqn 1.1.2}), together with the relation
$C_{n,k}=M_{2n-k, k}(0,0)$, where
\begin{eqnarray}\label{eqn 2.1.5}
\hskip1cm M_{n,k}(0,0) \hskip-.22cm &=&\hskip-.22cm
\left\{
\begin{array}{ll}
\frac{k+1}{n+1}\binom{n+1}{\frac{n-k}{2}}, & \mbox{if}\ n-k\ \mbox{ even }, \\[5pt]
0, & \mbox{otherwise}.
\end{array}\right. \hskip1cm\mbox{\cite[Example (iv)]{SunMa}}
\end{eqnarray}
we can get (\ref{eqn 2.1.1}). After some simple computation, one can easily derive (\ref{eqn 2.1.3}) from (\ref{eqn 2.1.1}).

Similarly, the case $(x,y)=(0,0)$ and $r=1$, after replacing $n,m$ respectively by $2n+1, 2m-\ell$ in (\ref{eqn 1.1.2}), reduces to (\ref{eqn 2.1.2}),
which, by some routine simplification, leads to (\ref{eqn 2.1.4}).\qed\vskip.2cm

When $m$ and $\ell$ in $\lambda_{n,k}(m,\ell)$ and $\mu_{n,k}(m,\ell)$ are assigned by special values, one can deduce that
\begin{eqnarray*}
\lambda_{n,k}(n,0)   \hskip-.22cm &=&\hskip-.22cm 2k(2n+1)(n+k+2),\\
\lambda_{n,k}(n,1)   \hskip-.22cm &=&\hskip-.22cm (n+k+2)(8nk+2k+2n+2),\\
\lambda_{n,k}(n+1,0) \hskip-.22cm &=&\hskip-.22cm (2k+3)(n-k+1)(2n+2), \\
\lambda_{n,k}(n+1,1) \hskip-.22cm &=&\hskip-.22cm (2k+2)(2n+1)(2n+2),\\
\lambda_{n,k}(n+2,1) \hskip-.22cm &=&\hskip-.22cm (n-k+1)(8nk+10k+14n+16);\\
\mu_{n,k}(n,0)       \hskip-.22cm &=&\hskip-.22cm (2k+1)(2n+2)(n+k+3), \\
\mu_{n,k}(n,1)   \hskip-.22cm &=&\hskip-.22cm (n+k+3)(8nk+6k+6n+6),\\
\mu_{n,k}(n+1,0)     \hskip-.22cm &=&\hskip-.22cm (2k+4)(n-k+1)(2n+3),\\
\mu_{n,k}(n+1,1) \hskip-.22cm &=&\hskip-.22cm (2k+3)(2n+2)(2n+3),\\
\mu_{n,k}(n+2,1) \hskip-.22cm &=&\hskip-.22cm (n-k+1)(8nk+14k+18n+30).
\end{eqnarray*}

In the cases $\lambda_{n,k}(m,1)$ for $m=n, n+1$ and $n+2$, after shifting $n$ to $n+1$ in the case $\lambda_{n,k}(n,1)$,
(\ref{eqn 2.1.3}) reproduces the known results \cite[Corollary 3.3]{SunMa}.
\begin{corollary} For any integer $n\geq 0$, there hold
\begin{eqnarray}
C_{n+1}^2 \hskip-.22cm &=&\hskip-.22cm
\sum_{k=0}^{n}\frac{ (2k+1)(2k+3)(8nk+10k+2n+4)}{(2n+1)_5}\binom{2n+2}{n-k}\binom{2n+5}{n-k+2},\nonumber\\
C_{n}C_{n+1} \hskip-.22cm &=&\hskip-.22cm
\sum_{k=0}^{n}\frac{(2k+1)(2k+2)(2k+3)}{(2n+1)(2n+2)(2n+3)^2}\binom{2n+3}{n-k}\binom{2n+3}{n-k+1}, \label{eqn 2.3.1}\\
C_{n}C_{n+2} \hskip-.22cm &=&\hskip-.22cm
\sum_{k=0}^{n}\frac{ (2k+1)(2k+3)(8nk+10k+14n+16)}{(2n+1)_5}\binom{2n+2}{n-k}\binom{2n+5}{n-k+1}.\nonumber
\end{eqnarray}
\end{corollary}
Note that another two cases $\lambda_{n,k}(n+1,0), \mu_{n,k}(n,0)$ also generate (\ref{eqn 2.3.1}).

In the cases $\mu_{n,k}(m,1)$ for $m=n, n+1$ and $n+2$, after shifting $n$ to $n+1$ in the case $\mu_{n,k}(n,1)$,
(\ref{eqn 2.1.4}) reproduces the known results \cite[Corollary 3.4]{SunMa}.
\begin{corollary} For any integer $n\geq 0$, there hold
\begin{eqnarray}
C_{n+1}C_{n+2}\hskip-.22cm &=&\hskip-.22cm
\sum_{k=0}^{n}\frac{(2k+2)(2k+4)(8nk+6n+14k+12)}{(2n+2)_5}\binom{2n+3}{n-k}\binom{2n+6}{n-k+2}, \nonumber\\
C_{n+1}^2 \hskip-.22cm &=&\hskip-.22cm
\sum_{k=0}^{n}\frac{(2k+2)(2k+3)(2k+4)}{(2n+2)(2n+3)(2n+4)^2}\binom{2n+4}{n-k}\binom{2n+4}{n-k+1}, \label{eqn 2.4.1} \\
C_{n+1}C_{n+2}\hskip-.22cm &=&\hskip-.22cm
\sum_{k=0}^{n}\frac{(2k+2)(2k+4)(8nk+14k+18n+30)}{(2n+2)_5}\binom{2n+3}{n-k}\binom{2n+6}{n-k+1}. \nonumber
\end{eqnarray}
\end{corollary}
Note that another two cases $\mu_{n,k}(n+1,0), \lambda_{n,k}(n,0)$ also generate (\ref{eqn 2.4.1}).

Taking $m=n$ in (\ref{eqn 2.1.3}) and (\ref{eqn 2.1.4}) into account, we have
\begin{eqnarray*}
\lambda_{n,k}(n,\ell)  \hskip-.22cm &=&\hskip-.22cm   (\ell+1)(n+k+2)(2k(2n+1)+\ell(n-k+1)), \\
\mu_{n,k}(n,\ell)  \hskip-.22cm &=&\hskip-.22cm   (\ell+1)(n+k+3)((2k+1)(2n+2)+\ell(n-k+1)),
\end{eqnarray*}
which yield the following results.
\begin{corollary} For any integers $n\geq\ell\geq 0$, there hold
\begin{eqnarray}
\frac{1}{2n-\ell+1}\binom{2n-\ell+1}{n-\ell}C_{n} \hskip-.22cm &=&\hskip-.22cm
\sum_{k=0}^{n-\ell}\frac{ (2k+1)(2k+\ell+2)\overline{\lambda}_{n,k}(\ell)}{(2n+1)_2(2n-\ell)_3}\binom{2n+2}{n-k+1}\binom{2n-\ell+2}{n-k-\ell},  \label{eqn 2.5.1}\\
\frac{1}{2n-\ell+1}\binom{2n-\ell+1}{n-\ell}C_{n+1} \hskip-.22cm &=&\hskip-.22cm
\sum_{k=0}^{n-\ell}\frac{ (2k+2)(2k+\ell+3)\overline{\mu}_{n,k}(\ell)}{(2n+2)_2(2n-\ell+1)_3}\binom{2n+3}{n-k+1}\binom{2n-\ell+3}{n-k-\ell}, \label{eqn 2.5.2}
\end{eqnarray}
where $\overline{\lambda}_{n,k}(\ell)=2k(2n+1)+\ell(n-k+1)$ and $\overline{\mu}_{n,k}(\ell)=(2k+1)(2n+2)+\ell(n-k+1)$.
\end{corollary}

It should be pointed out that both (\ref{eqn 2.5.1}) and (\ref{eqn 2.5.2}) are still
correct for any integer $\ell\leq -1$ if one notices that they
hold trivially for any integer $\ell> n$ and both sides of them can
be transferred into polynomials on $\ell$. Specially, after shifting $n$ to $n-1$,
the case $\ell=-1$ in (\ref{eqn 2.5.1}) and (\ref{eqn 2.5.2}) generates
the following corollary.
\begin{corollary} For any integer $n\geq 1$, there hold
\begin{eqnarray*}
\binom{2n}{n}C_{n-1}\hskip-.22cm &=&\hskip-.22cm
\sum_{k=0}^{n}\frac{(2k+1)^2(4nk-n-k)}{(2n-1)^2(2n)(2n+1)}\binom{2n}{n-k}\binom{2n+1}{n-k}, \\
\binom{2n}{n}C_{n}\hskip-.22cm &=&\hskip-.22cm
\sum_{k=0}^{n}\frac{(k+1)^2(4nk+n+k)}{n(n+1)(2n+1)^2}\binom{2n+1}{n-k}\binom{2n+2}{n-k}.
\end{eqnarray*}
\end{corollary}

\section{The third transformation on the Catalan triangle }

Let $\mathcal{Z}=(Z_{n,k})_{n\geq k\geq 0}$ be the infinite lower triangle defined on the Catalan triangle $\mathcal{C}$ by
\begin{eqnarray*}
Z_{2n,2k}=C_{n+k,2k}C_{n+k+1,2k+1},\ \ \  Z_{2n,2k+1}=C_{n+k+1,2k+1}C_{n+k+1,2k+2},\ \  (n\geq 0), \\
Z_{2n-1,2k}=C_{n+k,2k}C_{n+k,2k+1},\ \ \  Z_{2n-1,2k+1}=C_{n+k,2k+1}C_{n+k+1,2k+2},\ \  (n\geq 1).
\end{eqnarray*}
\begin{center}
\begin{eqnarray*}
\begin{array}{c|ccccccc|c|c}\hline
n/k & 0   & 1   & 2    & 3    & 4    & 5    & 6     & row\ sums & alternating\ sums\ of\ rows\   \\\hline
  0 & 1   &     &      &      &      &      &       &  1      &  1=1^2     \\
  1 & 1   & 1   &      &      &      &      &       &  2      &  0\hskip.8cm   \\
  2 & 2   & 2   & 1    &      &      &      &       &  5      &  1=1^2     \\
  3 & 4   & 6   & 3    & 1    &      &      &       &  14     &  0\hskip.8cm   \\
  4 & 10  & 15  & 12   & 4    &  1   &      &       &  42     &  4=2^2     \\
  5 & 25  & 45  & 36   & 20   &  5   &  1   &       &  132    &  0 \hskip.8cm   \\
  6 & 70  & 126 & 126  & 70   &  30  &  6   &  1    &  429    &  25=5^2    \\\hline
\end{array}
\end{eqnarray*}
Table 3.1. The values of $Z_{n,k}$ for $n$ and $k$ up to $6$, together with the row sums and the alternating sums of rows.
\end{center}
Table 3.1 illustrates the triangle $\mathcal{Z}$ for small $n$ and $k$ up to $6$, together
with the row sums and the alternating sums of rows. It signifies that
the sums and the alternating sums of rows of $\mathcal{Z}$ are in direct contact with the first column of $\mathcal{C}$.
Generally, we obtain the second result.

\begin{theorem} For any integers $m, n\geq 0$, let $p=m-n+1$. Then there hold
\begin{eqnarray}
\sum_{k=0}^{\min\{m,n\}}C_{m+k+1,2k+1}(C_{n+k,2k}+C_{n+k+1,2k+2})
\hskip-.22cm &=&\hskip-.22cm C_{m+n+1},  \label{eqn 3.1.1} \\
\sum_{k=0}^{\min\{m,n\}} C_{m+k+1,2k+1}(C_{n+k,2k}-C_{n+k+1,2k+2})
\hskip-.22cm &=&\hskip-.22cm G_{n,m}(p), \label{eqn 3.1.2}
\end{eqnarray}
where
\begin{eqnarray*}
G_{n,m}(p)=\left\{\begin{array}{rl}
\sum_{i=0}^{p-1}C_{n+i}C_{m-i}, & if  \ p\geq 1, \\
0,                              & if  \ p=0,      \\
-\sum_{i=1}^{|p|}C_{n-i}C_{m+i}, & if \ p\leq -1,
\end{array}\right.
\end{eqnarray*}
\end{theorem}
\pf (\ref{eqn 3.1.1}) is equivalent to (\ref{eqn 1.1.1}), if one notices that
$$C_{m+k+1,2k+1}=B_{m,k}, B_{n,k}=C_{n+k+1,2k+1}=C_{n+k,2k}+C_{n+k+1,2k+2}.$$

For the case $p=m-n+1\geq 0$ in (\ref{eqn 3.1.2}),
setting $(x,y)=(0,0)$ in (\ref{eqn 1.1.2}), together with the relation $C_{n,k}=M_{2n-k, k}(0,0)$ and (\ref{eqn 2.1.5}),
we have
\begin{eqnarray*}
\lefteqn{\sum_{k=0}^{\min\{m,n\}} C_{m+k+1,2k+1}(C_{n+k,2k}-C_{n+k+1,2k+2}) }\\
\hskip-.22cm &=&\hskip-.22cm
\sum_{k=0}^{n}\det\left(\begin{array}{cc}
C_{n+k,2k}   & 0 \\[5pt]
0            & C_{n+p+k,2k+1}
\end{array}\right)+ \det\left(\begin{array}{cc}
            0   & C_{n+k+1,2k+2} \\[5pt]
C_{n+p+k,2k+1}  & 0
\end{array}\right) \\
\hskip-.22cm &=&\hskip-.22cm
\sum_{k=0}^{2n}\det\left(\begin{array}{cc}
M_{2n,k}(0,0)        & M_{2n,k+1}(0,0) \\[5pt]
M_{2n+2p-1,k}(0,0) & M_{2n+2p-1,k+1}(0,0)
\end{array}\right) \\
\hskip-.22cm &=&\hskip-.22cm \sum_{i=0}^{2p-2} M_{2n+i,0}(0,0)M_{2n+2p-i-2,0}(0,0)=\sum_{i=0}^{p-1} M_{2n+2i,0}(0,0)M_{2n+2p-2i-2,0}(0,0) \\
\hskip-.22cm &=&\hskip-.22cm \sum_{i=0}^{p-1}C_{n+i}C_{n+p-i-1}=\sum_{i=0}^{p-1}C_{n+i}C_{m-i}, \\
\end{eqnarray*}
as desired.

Similarly, the case $p\leq -1$ can be proved, the details are left to interested readers.\qed\vskip.2cm

Note that a weighted partial Motzkin path with no horizontal steps is just a partial Dyck path. Then the relation $C_{n,k}=M_{2n-k,k}(0,0)$ signifies
that $C_{n,k}$ counts the set $\mathcal{C}_{n,k}$ of partial Dyck paths of length $2n-k$ from $(0,0)$ to $(2n-k, k)$ \cite{MerliSprug}.
Such partial Dyck paths have exactly $n$ up steps and $n-k$ down steps. For any step, we say that it is at level
$i$ if the $y$-coordinate of its end point is $i$.
If $P=L_1L_2\dots L_{2n-k-1-1}L_{2n-k}\in \mathcal{C}_{n,k}$, denote by
$\overline{P}=\overline{L}_{2n-k}\overline{L}_{2n-k-1}\dots
\overline{L}_2\overline{L}_1$ the reverse of the path $P$, where $\overline{L}_i=\mathbf{u}$ if $L_i=\mathbf{d}$ and $\overline{L}_i=\mathbf{d}$
if $L_i=\mathbf{u}$.

For $k=0$, a partial Dyck path  is an (ordinary) Dyck path. For any Dyck path $P$ of length $2n+2m+2$,
its $(2n+1)$-th step $L$ (along the path) must end at odd level, say $2k+1$ for some $k\geq 0$, then $P$ can be uniquely partitioned
into $P=P_1LP_2$, where $(P_1, \overline{P_2})\in \mathcal{C}_{n+k,2k}\times\mathcal{C}_{m+k+1,2k+1}$
if $L=\mathbf{u}$ and $(P_1, \overline{P_2})\in \mathcal{C}_{n+k+1,2k+2}\times\mathcal{C}_{m+k+1,2k+1}$
if $L=\mathbf{d}$. Hence, the cases $p=0, 1$ and $2$, i.e., $m=n-1, n$ and $n+1$ in (\ref{eqn 3.1.2}) produce the following corollary.

\begin{corollary} For any integer $n\geq 0$, according to the $(2n+1)$-th step $\mathbf{u}$ or $\mathbf{d}$, we have
\begin{itemize}
\item[(i)] The number of Dyck paths of length $4n$ is bisected;
\item[(ii)] The parity of Dyck paths of length $4n+2$ is $C_n^2$;
\item[(iii)] The parity of Dyck paths of length $4n+4$ is $2C_nC_{n+1}$.
\end{itemize}
\end{corollary}

\section{The forth transformation on the Catalan triangle }

Let $\mathcal{W}=(W_{n,k})_{n\geq k\geq 0}$ be the infinite lower triangle defined on the Catalan triangle by
\begin{eqnarray*}
W_{2n,k} \hskip-.22cm &=&\hskip-.22cm \rm{per}\left(\begin{array}{cc}
C_{n+k, 2k}    & C_{n+k,2k+1}  \\[5pt]
C_{n+k+1, 2k}  & C_{n+k+1,2k+1}
\end{array}\right), \\
W_{2n+1,k} \hskip-.22cm &=&\hskip-.22cm {\rm per}\left(\begin{array}{cc}
C_{n+k, 2k}    & C_{n+k,2k+1}  \\[5pt]
C_{n+k+2, 2k}  & C_{n+k+2,2k+1}
\end{array}\right),
\end{eqnarray*}
where ${\rm per}(A)$ denotes the permanent of a square matrix $A$. Table 4.1 illustrates the triangle $\mathcal{W}$
for small $n$ and $k$ up to $8$, together with the row sums.

\begin{center}
\begin{eqnarray*}
\begin{array}{c|ccccccccc|cc}\hline
n/k & 0     & 1     & 2      & 3      & 4    & 5    & 6  & 7    & 8    & row\ sums        \\\hline
  0 & 1     &       &        &        &      &      &    &      &      &  1           \\
  1 & 2     &       &        &        &      &      &    &      &      &  2           \\
  2 & 4     & 1     &        &        &      &      &    &      &      &  5           \\
  3 & 10    & 4     &        &        &      &      &    &      &      &  14            \\
  4 & 20    & 21    & 1      &        &      &      &    &      &      &  42            \\
  5 & 56    & 70    & 6      &        &      &      &    &      &      &  132           \\
  6 & 140   & 238   & 50     & 1      &      &      &    &      &      &  429            \\
  7 & 420   & 792   & 210    & 8      &      &      &    &      &      &  1430           \\
  8 & 1176  & 2604  & 990    & 91     &  1   &      &    &      &      &  4862          \\\hline
\end{array}
\end{eqnarray*}
Table 4.1. The values of $W_{n,k}$ for $n$ and $k$ up to $8$, together with the row sums.
\end{center}

This motivates us to consider the permanent transformation on the triangle $\mathcal{M}=(M_{n,k}(x,y))_{n\geq k\geq 0}$.
Recall that $M_{n,k}(x,y)$ is the weight sum of the set $\mathcal{M}_{n,k}(x,y)$ of
all weighted partial Motzkin paths ending at $(n,k)$. For any step of a partial weighted Motzkin path $P$, we say that it is at level
$i$ if the $y$-coordinate of its end point is $i$.
If $P=L_1L_2\dots L_{2n-k-1-1}L_{2n-k}\in \mathcal{M}_{n,k}(x,y)$, denote by
$\overline{P}=\overline{L}_{2n-k}\overline{L}_{2n-k-1}\dots
\overline{L}_2\overline{L}_1$ the reverse of the path $P$, where $\overline{L}_i=\mathbf{u}, \mathbf{h}$ or $\mathbf{d}$
if $L_i=\mathbf{d}, \mathbf{h}$ or $\mathbf{u}$ respectively.
An up step $\mathbf{u}$ of $P$ at level $i$ is $R$-$visible$ if it is the rightmost
up step at level $i$ and there are no other up steps at the same
level to its right.

\begin{theorem}
For any integers $m, n, r$ with $m\geq n\geq 0$, there
holds
\begin{eqnarray}\label{eqn 4.1.1}
\sum_{k=0}^{m}\rm{per}\left(\begin{array}{cc}
M_{n,k}(y,y) & M_{n+r,k+1}(y,y) \\[5pt]
M_{m,k}(y,y) & M_{m+r,k+1}(y,y)
\end{array}\right)
\hskip-.22cm &=&\hskip-.22cm M_{m+n+r,1}(y,y)+H_{n,m}(r),
\end{eqnarray}
where
\begin{eqnarray*}
H_{n,m}(r)=\left\{\begin{array}{rl}
\sum_{i=0}^{r-1}M_{n+i,0}(y,y)M_{m+r-i-1,0}(y,y),    & if \ r\geq 1, \\
0,                                                   & if \ r=0,      \\
-\sum_{i=1}^{|r|}M_{n-i,0}(y,y)M_{m-|r|+i-1,0}(y,y), & if \ r\leq
-1.
\end{array}\right.
\end{eqnarray*}
\end{theorem}

\pf We just give the proof of the part when $r\geq 0$, the other
part can be done similarly and is left to interested readers. Define
\begin{eqnarray*}
{\mathcal{A}}_{n,m,k}^{(r)}  \hskip-.22cm &=&\hskip-.22cm \{(P, Q)|P\in \mathcal{M}_{n, k}(y,y), Q\in \mathcal{M}_{m+r, k+1}(y,y)\}, \\
{\mathcal{B}}_{n,m,k}^{(r)}  \hskip-.22cm &=&\hskip-.22cm \{(P,
Q)|P\in \mathcal{M}_{n+r, k+1}(y,y), Q\in \mathcal{M}_{m, k}(y,y)\},
\end{eqnarray*}
and ${\mathcal{C}}_{n,m,k}^{(r,i)}$ to be the subset of
${\mathcal{A}}_{n,m,k}^{(r)}$ such that for any $(P, Q)\in
{\mathcal{C}}_{n,m,k}^{(r,i)}$, $Q=Q_{1}\mathbf{u}Q_{2}$ with $Q_{1}\in
\mathcal{M}_{i,k}(y,y)$ and $Q_2\in \mathcal{M}_{m+r-i-1,0}(y,y)$
for $k\leq i\leq r-1$. In other words, for any $(P, Q)\in
{\mathcal{C}}_{n,m,k}^{(r,i)}$, $Q$ satisfies the conditions
that (a) the last R-visible up step of $Q$ stays at level $k+1$, and
(b) there are exactly $i$ steps immediately ahead of the last
R-visible up step of $Q$ for $k\leq i\leq r-1$. Clearly,
${\mathcal{C}}_{n,m,k}^{(r,i)}$ is the empty set if $r=0$.

It is easily to see that the weights of the sets
${\mathcal{A}}_{n,m,k}^{(r)}$ and
${\mathcal{B}}_{n,m,k}^{(r)}$ are
\begin{eqnarray*}
w({\mathcal{A}}_{n,m,k}^{(r)}) \hskip-.22cm &=&\hskip-.22cm M_{n, k}(y,y)M_{m+r, k+1}(y,y), \\
w({\mathcal{B}}_{n,m,k}^{(r)}) \hskip-.22cm &=&\hskip-.22cm
M_{n+r, k+1}(y,y)M_{m, k}(y,y).
\end{eqnarray*}
For $0\leq i\leq r-1$, the weight of the set
$\bigcup_{k=0}^{i}{\mathcal{C}}_{n,m,k}^{(r,i)}$ is
$M_{n+i,0}(y,y)M_{m+r-i-1,0}(y,y)$. This claim can be verified by
the following argument. For any $(P, Q)\in
{\mathcal{C}}_{n,m,k}^{(r,i)}$, we have $Q=Q_{1}\mathbf{u}Q_{2}$ as
mentioned above with $Q_{1}\in \mathcal{M}_{i,k}(y,y)$ and $Q_2\in
\mathcal{M}_{m+r-i-1,0}(y,y)$, then $P\overline{Q}_{1}\in
\mathcal{M}_{n+i,0}(y,y)$ such that the last $(i+1)$-th step of
$P\overline{Q}_{1}$ is at level $k$. Summing $k$ for $0\leq k\leq
i$, all $P\overline{Q}_{1}\in \mathcal{M}_{n+i,0}(y,y)$ contribute
the total weight $M_{n+i,0}(y,y)$ and all $Q_2\in \mathcal{M}_{m+r-i-1,0}(y,y)$ contribute the total weight
$M_{m+r-i-1,0}(y,y)$. Hence,
$w(\bigcup_{k=0}^{i}{\mathcal{C}}_{n,m,k}^{(r,i)})=M_{n+i,0}(y,y)M_{m+r-i-1,0}(y,y)$,
and then
\begin{eqnarray*}
w(\bigcup_{i=0}^{r-1}\bigcup_{k=0}^{i}{\mathcal{C}}_{n,m,k}^{(r,i)})
=w(\bigcup_{k=0}^{r-1}\bigcup_{i=k}^{r-1}{\mathcal{C}}_{n,m,k}^{(r,i)})=\sum_{i=0}^{r-1}
M_{n+i,0}(y,y)M_{m+r-i-1,0}(y,y)=H_{n,m}(r).
\end{eqnarray*}

Let
${\mathcal{A}}_{n,m}^{(r)}=\bigcup_{k=0}^{n+r-1}{\mathcal{A}}_{n,m,k}^{(r)}$,
${\mathcal{B}}_{n,m}^{(r)}=\bigcup_{k=0}^{m+r-1}{\mathcal{B}}_{n,m,k}^{(r)}$
and
${\mathcal{C}}_{n,m,k}^{(r)}=\bigcup_{i=k}^{r-1}{\mathcal{C}}_{n,m,k}^{(r,i)}$.
To prove (\ref{eqn 4.1.1}), it suffices to construct a bijection
$\varphi$ between
${\mathcal{B}}_{n,m}^{(r)}\bigcup\big({\mathcal{A}}_{n,m}^{(r)}-
\bigcup_{k=0}^{r-1}{\mathcal{C}}_{n,m,k}^{(r)}\big)$ and
$\mathcal{M}_{m+n+r,1}(y,y)$.

For any $(P,Q)\in {\mathcal{B}}_{n,m,k}^{(r)}$, $P\overline{Q}$
is exactly an element of $\mathcal{M}_{m+n+r,1}(y,y)$. Note that in
this case, the first R-visible up step of $P$ is still the one of
$P\overline{Q}$ and it is at most the $(n+r)$-th step of
$P\overline{Q}$.

For any $(P,Q)\in
{\mathcal{A}}_{n,m,k}-{\mathcal{C}}_{n,m,k}^{(r)}$, find the
last R-visible up step $\mathbf{u}^{*}$ of $Q$, $Q$ can be uniquely
partitioned into $Q=Q_1\mathbf{u}^{*}Q_2$, where $Q_1\in
\mathcal{M}_{j,k}(y,y)$ for some $j\geq r$, then
$P\overline{Q}_1\mathbf{u}^{*} Q_2$ forms an element of
$\mathcal{M}_{m+n+r,1}(y,y)$. Note that in this case, the last
R-visible up step $\mathbf{u}^{*}$ of $Q$ is still the one of
$P\overline{Q}_1\mathbf{u}^{*}Q_2$. Moreover, the $\mathbf{u}^{*}$ step is at least
the $(n+r+1)$-th step of $P\overline{Q}_1\mathbf{u}^{*}Q_2$.

Conversely, for any path in $\mathcal{M}_{m+n+r,1}(y,y)$, it can be
partitioned uniquely into $PQ$, where $P\in \mathcal{M}_{n+r,
k}(y,y)$ for some $k\geq 0$. If the unique R-visible up step $\mathbf{u}^{*}$
of $PQ$ is lying in $P$, then $k\geq 1$ and $(P, \overline{Q})\in
{\mathcal{B}}_{n,m,k-1}^{(r)}$; If the $\mathbf{u}^{*}$ step is lying in
$Q$, $PQ$ can be repartitioned into $P_1P_2\mathbf{u}^{*}Q_1$ with $P_1\in
\mathcal{M}_{n,j}(y,y)$ for some $j\geq 0$, then $(P_1,
\overline{P_2}\mathbf{u}^{*}Q_1)\in
{\mathcal{A}}_{n,m,j}-{\mathcal{C}}_{n,m,j}^{(r)}$.

Clearly, the above procedure is invertible. Hence, $\varphi$ is
indeed a bijection as desired and (\ref{eqn 4.1.1}) is proved. \qed\vskip0.2cm

\begin{theorem} For any integers $m, n, p$ with $m\geq n\geq 0$, there hold
\begin{eqnarray}
\sum_{k=0}^{m}\rm{per}\left(\begin{array}{cc}
C_{n+k, 2k}  & C_{n+p+k,2k+1}  \\[5pt]
C_{m+k, 2k}  & C_{m+p+k,2k+1}
\end{array}\right)
\hskip-.22cm &=&\hskip-.22cm C_{m+n+p,1}+F_{n,m}(p), \label{eqn 4.1.2} \\
\sum_{k=0}^{m}\rm{per}\left(\begin{array}{cc}
C_{n+k, 2k+1}  & C_{n+p+k+1,2k+2}  \\[5pt]
C_{m+k, 2k+1}  & C_{m+p+k+1,2k+2}
\end{array}\right)
\hskip-.22cm &=&\hskip-.22cm C_{m+n+p,1}+F_{n,m}(p), \label{eqn 4.1.3}
\end{eqnarray}
where
\begin{eqnarray*}
F_{n,m}(p)=\left\{\begin{array}{rl}
\sum_{i=0}^{p-1}C_{n+i}C_{m+p-i-1}, & if \ p\geq 1, \\
0,                                  & if \ p=0,      \\
-\sum_{i=1}^{|p|}C_{n-i}C_{m-|p|+i-1}, & if \ p\leq -1.
\end{array}\right.
\end{eqnarray*}
\end{theorem}
\pf To prove (\ref{eqn 4.1.2}), replacing $n, m, r$ respectively by $2n, 2m, 2p-1$ and
setting $(y,y)=(0,0)$ in (\ref{eqn 4.1.1}), together with the relation $C_{n,k}=M_{2n-k, k}(0,0)$ and (\ref{eqn 2.1.5}),
we have
\begin{eqnarray*}
\lefteqn{\sum_{k=0}^{m}\rm{per}\left(\begin{array}{cc}
C_{n+k, 2k}  & C_{n+p+k,2k+1}  \\[5pt]
C_{m+k, 2k}  & C_{m+p+k,2k+1}
\end{array}\right)}\\
\hskip-.22cm &=&\hskip-.22cm
\sum_{k=0}^{m}\rm{per}\left(\begin{array}{cc}
M_{2n,2k}(0,0)   & M_{2n+2p-1, 2k+1}(0,0) \\[5pt]
M_{2m,2k}(0,0)   & M_{2m+2p-1, 2k+1}(0,0)
\end{array}\right)\\
\hskip-.22cm &=&\hskip-.22cm
\sum_{k=0}^{2m}\rm{per}\left(\begin{array}{cc}
M_{2n,k}(0,0)   & M_{2n+2p-1, k+1}(0,0) \\[5pt]
M_{2m,k}(0,0)   & M_{2m+2p-1, k+1}(0,0)
\end{array}\right) \\
\hskip-.22cm &=&\hskip-.22cm M_{2n+2m+2p-1,1}(0,0)+H_{2n,2m}(2p-1)\\
\hskip-.22cm &=&\hskip-.22cm C_{m+n+p,1}+F_{n,m}(p),
\end{eqnarray*}
as desired.

Similarly, replacing $n, m, r$ respectively by $2n-1, 2m-1, 2p+1$ and
setting $(y,y)=(0,0)$ in (\ref{eqn 4.1.1}), together with the relation $C_{n,k}=M_{2n-k, k}(0,0)$
and (\ref{eqn 2.1.5}), one can prove (\ref{eqn 4.1.3}), the details are left to interested readers.\qed\vskip.2cm

The case $p=0$ in (\ref{eqn 4.1.2}) and (\ref{eqn 4.1.3}), after some routine computation, generates
\begin{corollary} For any integers $m\geq n\geq 1$, there hold
\begin{eqnarray*}
C_{n+m} \hskip-.22cm &=&\hskip-.22cm
\sum_{k=0}^{n} \frac{(2k+1)(2k+2)(4mn-2(m+n)k)}{(2n)(2n+1)(2m)(2m+1)
}\binom{2n+1}{n-k}\binom{2m+1}{m-k}, \\
C_{n+m} \hskip-.22cm &=&\hskip-.22cm
\sum_{k=0}^{n-1} \frac{(2k+2)(2k+3)(4mn+4m+4n+2(m+n)k)}{(2n)(2n+1)(2m)(2m+1)
}\binom{2n+1}{n-k-1}\binom{2m+1}{m-k-1}.
\end{eqnarray*}
Specially, the $m=n$ case produces
\begin{eqnarray*}
C_{2n} \hskip-.22cm &=&\hskip-.22cm \sum_{k=0}^{n-1} \frac{(2k+1)(2k+2)}{n(2n+1)}\binom{2n}{n-k-1}\binom{2n+1}{n-k}, \\
C_{2n} \hskip-.22cm &=&\hskip-.22cm \sum_{k=0}^{n-1} \frac{(2k+2)(2k+3)}{n(2n+1)}\binom{2n}{n-k-1}\binom{2n+1}{n-k-1}.
\end{eqnarray*}
\end{corollary}

The case $p=1$ in (\ref{eqn 4.1.2}), after some routine computation, generates
\begin{corollary} For any integers $m\geq n\geq 0$, there hold
\begin{eqnarray*}
C_{n+m+1}+C_{n}C_{m} \hskip-.22cm &=&\hskip-.22cm
\sum_{k=0}^{n} \frac{(2k+1)(2k+2)\eta_{n,m}(k)}{(2n+1)(2n+2)(2m+1)(2m+2)
}\binom{2n+2}{n-k}\binom{2m+2}{m-k},
\end{eqnarray*}
where $\eta_{n,m}(k)=4mn+5(m+n)+2(m+n+1)k+4$. Specially, the $m=n$ case produces
\begin{eqnarray*}
C_{2n+1}+C_{n}^2 \hskip-.22cm &=&\hskip-.22cm \sum_{k=0}^{n} \frac{(2k+1)(2k+2)}{(n+1)(2n+1)}\binom{2n+1}{n-k}\binom{2n+2}{n-k}.
\end{eqnarray*}
\end{corollary}

In the case $y=2$ and $r=p$ in (\ref{eqn 4.1.1}), together with the relations $B_{n,k}=M_{n,k}(2,2)$ and $B_{n,0}=C_{n+1}$, similar to
the proof of (\ref{eqn 4.1.2}), we obtain a result on Shapiro's Catalan triangle.
\begin{theorem} For any integers $m, n, p$ with $m\geq n\geq 0$, there holds
\begin{eqnarray}\label{eqn 4.1.4}
\sum_{k=0}^{m}\rm{per}\left(\begin{array}{cc}
B_{n,k}  & B_{n+p,k+1}    \\[5pt]
B_{m,k}  & B_{m+p,k+1}
\end{array}\right)
\hskip-.22cm &=&\hskip-.22cm B_{m+n+p,1}+F_{n+1,m+1}(p).
\end{eqnarray}
\end{theorem}

The case $p=0$ in (\ref{eqn 4.1.4}), after some routine computation, generates
\begin{corollary} For any integers $m\geq n\geq 0$, there holds
\begin{eqnarray*}
\frac{2}{n+m+1}\binom{2n+2m+2}{n+m-1} \hskip-.22cm &=&\hskip-.22cm
\sum_{k=0}^{n} \frac{(2k+2)(2k+4)\nu_{n,k}(m)}{(2n+2)_2(2m+2)_2
}\binom{2n+3}{n-k}\binom{2m+3}{m-k}.
\end{eqnarray*}
where $\nu_{n,k}(m)=2mn+3m+3n-6k-2k^2$. Specially, the $m=n$ case
produces
\begin{eqnarray*}
\frac{1}{2n+1}\binom{4n+2}{2n-1} \hskip-.22cm &=&\hskip-.22cm
\sum_{k=0}^{n-1} \frac{(k+1)(k+2)}{(n+1)^2
}\binom{2n+2}{n-k-1}\binom{2n+2}{n-k}.
\end{eqnarray*}
\end{corollary}

\vskip1cm
\section*{Acknowledgements} {  The work
was partially supported by The National Science Foundation of China
and by the Fundamental Research Funds for the Central Universities.}

%==============================================================================================================


\begin{thebibliography}{99}



\bibitem{AignerA} M. Aigner, {\it Catalan-like numbers and determinants}, J. Combin. Theory Ser.
A, 87 (1999), 33-51.

\bibitem{AignerB} M. Aigner, {\it Catalan and other numbers -- a recurrent theme},
in: H. Crapo, D. Senato (Eds.), Algebraic Combinatorics and Computer Science, Springer,
Berlin (2001), 347-390.

\bibitem{AignerC} M. Aigner, {\it Enumeration via ballot numbers}, Discrete Math., 308 (2008), 2544-2563.

\bibitem{BacchMerli} D. Baccherini, D. Merlini, and R. Sprugnoli, {\it Level generating trees and proper Riordan arrays}, Applicable
Analysis and Discrete Mathematics, 2 (2008), 69-91.

%\bibitem{BarcPinz} E. Barcucci, R. Pinzani and R. Sprugnoli, {\it The Motzkin family}, Pure Math. Appl. 2 (1991), 249-279.

\bibitem{BarcVerri} E. Barcucci and M.C. Verri, {\it Some more properties of Catalan numbers},
Discrete Math., 102(3) (1992), 229-237.



%\bibitem{CamNk} N. Cameron and A. Nkwanta, {\it On some (pseudo) involutions in the Riordan
%group}, J. Integer Seq., 8 (2005), Article 05.3.7.

\bibitem{ChenChu} X. Chen and W. Chu, {\it Moments on Catalan number}, J. Math. Anal. Appl., 349 (2) (2009), 311-316.

\bibitem{ChenLi} W.Y.C. Chen, N.Y. Li, L.W. Shapiro and S.H.F. Yan,
{\it Matrix identities on weighted partial Motzkin paths}, Europ. J.
Combin., 28 (2007), 1196-1207.

\bibitem{CheonJin} G.-S. Cheon, S.-T. Jin, {\it Structural properties of Riordan matrices and extending the
matrices}, Linear Algebra and its Appl., 435 (2011), 2019-2032.

\bibitem{CheonKim} G.-S. Cheon and H. Kim, {\it Simple proofs of open problems about the structure
of involutions in the Riordan group}, Linear Algebra and its Appl., 428 (2008), 930-940.

\bibitem{CheonKimShap} G.-S. Cheon, H. Kim and L.W. Shapiro, {\it Combinatorics of Riordan arrays
with identical A and Z sequences}, Discrete Math., 312(12-13) (2012), 2040-2049.

%\bibitem{Comtet} L. Comtet, Advanced Combinatorics, D. Reidel, Dordrecht, 1974.

%\bibitem{Deng} E.Y.P. Deng and W.-J. Yan, {\it Some identities on the Catalan, Motzkin and Schr$\ddot{o}$der
%numbers}, Discrete Applied Mathematics, 156(14), 28 (2008), 2781-2789.

%\bibitem{Deutsch} E. Deutsch, {\it Dyck path enumeration}, Discrete Math., 204 (1999), 167-202.


\bibitem{Eplett} W.J.R. Eplett, {\it A note about the Catalan triangle}, Discrete Math., 25 (1979), 289-291.


\bibitem{EuLiuYeh} S.-P. Eu, S.-C. Liu and Y.-N. Yeh, {\it Taylor expansions for Catalan and Motzkin
numbers}, Adv. Applied Math., 29 (2002), 345-357.

\bibitem{Forder} H.G. Forder, {\it Some problems in combinatorics}, Math. Gazette, 45 (1961), 199-201.

\bibitem{Gessel} I. Gessel, {\it Super ballot numbers}, J. Symbolic Comput., 14 (1992), 179-194.

\bibitem{GuoZeng} V.J.W. Guo and J. Zeng, {\it Factors of binomial sums from Catalan triangle}, J. Number Theory, 130 (1) (2010), 172-186.

\bibitem{Hilton} P. Hilton and J. Pedersen, {\it Catalan numbers, their generalization and their uses}, Math. Intelligencer, 13 (1991), 64-75.

\bibitem{Gutierrez} J.M. Gutierrez, M.A. Hern¨¢ndez, P.J. Miana, N. Romero, {\it New identities in the Catalan triangle},
J. Math. Anal. Appl. 341 (1) (2008) 52-61.


\bibitem{KitLies} S. Kitaev and J. Liese, {\it Harmonic numbers, Catalan's triangle and mesh patterns},
Discrete Math., 313(14) (2013), 1515-1531.

\bibitem{Knuth} D.E. Knuth, The Art of Computer Programming, 3rd edn., Addison-Wesley, 1998.

\bibitem{Luzona} A. Luzona, D. Merlini, M.A. Moronc and R. Sprugnoli, {\it Identities induced by Riordan arrays},
Linear Algebra and its Appl., 436 (2012), 631-647.


\bibitem{ShiMeiMa} S.-M. Ma, {\it Some combinatorial arrays generated by context-free grammars},
Europ. J. Combinatorics, 34(7) (2013), 1081-1091.

\bibitem{Merlini} D. Merlini, {\it Proper generating trees and their internal path length}, Discrete Applied
Math., 156 (2008), 627-646.

\bibitem{MerliSprug} D. Merlini, R. Sprugnoli and M.C. Verri, {\it Some statistics on Dyck paths},
J. Statistical Planning and Inference, 101 (2002), 211-227.


\bibitem{MerlVer} D. Merlini, M. C.Verri, {\it Generating trees and proper Riordan arrays}, Discrete Math., 218 (2000), 167-183.

\bibitem{Miana} P.J. Miana and N. Romero, {\it Computer proofs of new identities in the Catalan triangle}, Biblioteca de la Revista Matem¨¢tica
Iberoamericana, in: Proceedings of the "Segundas Jornadas de Teor¨ªa de N¨²meros", (2007), 1-7.

\bibitem{MianaRom} P.J. Miana and N. Romero, {\it Moments of combinatorial and Catalan numbers}, J. Number Theory, 130 (2010), 1876-1887.

\bibitem{Petkov} M. Petkov$\breve{s}$ek, H.S. Wilf and D. Zeilberger, {\it A=B}, A. K. Peters, Wellesley, MA, 1996.

\bibitem{RogersA} D.G. Rogers, {\it Eplett's identities for renewal arrays}, Discrete Math., 36 (1981), 97-l02.

\bibitem{RogersB} D.G. Rogers, {\it Pascal triangles, Catalan numbers and renewal arrays}, Discrete Math., 22 (1978), 301-310.

%\bibitem{Sedg}  R. Sedgewick and P. Flajolet, {\it An Introduction to the Analysis of Algorithms}, Addison-Wesley, Reading, MA, 1996.

\bibitem{ShapA} L.W. Shapiro, {\it A Catalan triangle}, Discrete Math., 14 (1976), 83-90.

\bibitem{ShapB} L.W. Shapiro, {\it Bijections and the Riordan group}, Theoret. Comput. Sci., 307
(2003), 403-413.

\bibitem{ShapGet} L.W. Shapiro, S. Getu, W.-J. Woan, L.C. Woodson, {\it The Riordan group}, Discrete Appl. Math., 34 (1991), 229-239.

\bibitem{Sloane} N.J.A. Sloane, {\it The On-Line Encyclopedia of Integer Sequences}, \\
http://www.research.att.com/{$\sim$}njas/sequences.

\bibitem{Sprugnoli} R. Sprugnoli, {\it Combinatorial sums through Riordan arrays}, J. Geom. 101 (2011), 195-210.

\bibitem{Sprug} R. Sprugnoli, {\it Riordan arrays and combinatorial sums}, Discrete Math., 132
(1994), 267-290.

%\bibitem{StanleyA} R.P. Stanley, {\it Catalan Addendum}, http://www-math.mit.edu/~rstan/ec/
%catadd.pdf, version of 13, July 2012.

%\bibitem{StanleyB} R.P. Stanley, Enumerative Combinatorics, vol. 2, Cambridge University Press, Cambridge, New York, 1999.

\bibitem{SunMa} Y. Sun and L. Ma, {\it Minors of a class of Riordan arrays related to weighted partial Motzkin paths}, submitted.

\bibitem{Zeilberger} D. Zeilberger, {\it The method of creative telescoping}, J. Symbolic Comput.,
11 (1991), 195-204.

\bibitem{ZhangPang} Z. Zhang and B. Pang, {\it Several identities in the Catalan triangle}, Indian J. Pure Appl. Math., 41(2)(2010), 363-378.





\end{thebibliography}
\end{document}